\documentclass[a4paper,11pt]{article}
\usepackage{latexsym}
\usepackage{amsfonts,amssymb}
\usepackage{amscd}
\usepackage{amsmath}
\usepackage[dvips]{graphicx}

\begin{document}

\begin{center}

{\sc More on Algebraic Structure of the Complete Partition Function \\for the $
Z_n $ - Potts Model, Part 1}

\medskip
A. K. Kwa\'sniewski$^{*}$, W. Bajguz$^{**}$

\smallskip

\noindent {\small (*)the member of the Institute of Combinatorics and its Applications }\\
\noindent {\small  Faculty of Physics, Bia{\l }ystok University}\\
\noindent {\small ul. Lipowa 41,  15 424  Bia\l ystok, Poland}\\
\noindent {\small e-mail: kwandr@gmail.com}\\
{\small $^{**}$Bia{\l }ystok University,\\Institute of Computer Science, ul.
Sosnowa 64\\PL -- 15-887 Bia{\l }ystokk, Poland}

\end{center}

\subsubsection*{Abstract}

{\small In this first part of  a larger review undertaking the results of the first author and a part
of the second author  doctor dissertation are presented. Next we plan to give a survey of a nowadays
situation in the area of investigation. Here we report on what follows. 
 
Calculation of the partition function for any vector Potts model is
at first reduced to the calculation of  traces of products of the generalized Clifford algebra generators.
The formula for such traces is derived. }

{\small The latter enables one  , in principle, to use an explicit calculation algorithm for  partition
functions also in other models for which the transfer matrix is an element
from generalized Clifford algebra. 
 }

{\small \noindent The method - simple for $Z_2$ case - becomes complicated
for $Z_n$, $n>2$, however everything is controlled due to knowledge of the
corresponding algebra properties and those of generalized $\cosh $ function. }

{\small \noindent Hence the work to gain the thermodynamics of the system,
though possibly tantalous, looks now a reasonable, tangible task with help
of computer symbolic calculations. }

{\small \noindent The discussion of the content of the in statu nascendi second part
is to be found at the end of this Part 1 presentation. This constitutes the section V.}

\smallskip\
PACS numbers: 05.50.+q,

\smallskip\

Keywords:  Potts models, generalized Clifford Algebras

\subsection*{I. Introduction}

The main idea of all calculations to follow (see [8]) is to consider the
task of determining the complete partition function for {\bf noncritical}
Potts models - as a problem from the theory of generalized Clifford algebras 
$C_{2p}^{\left( n\right) }$ classified in [15].

\noindent As a matter of fact, the transfer matrix approach has led the
authors of [1,14 except for 20,21] to these very algebras, although this
observation does not seem to be realized by the mentioned authors.

\noindent In general - the transfer matrix technique for a statistical
system (or - a lattice field theory) with the most general translational
invariant and globally symmetric Hamiltonian (Action) on a two-dimensional
lattice - does generate appropriate algebras of operators which are of the
type of algebra extensions of some groups [8,18,24].

\noindent If one considers $Z_n$ cyclic groups as symmetry groups of
Hamiltonian (chiral Action) then the algebra generated by transfer matrix
approach is the corresponding $C_k^{\left( n\right) }$ algebra, and it is
due to its properties, that the given model has the duality property [18].

\noindent For duality property we refer the reader to the review [19], and
for Potts models, in general - to [2,23].

\noindent Calculations to be carried out here for the $Z_n$ - vector Potts
model simplify tremendously in the case of $n=2$ i.e. for the Ising model
and there lead to the known complete partition function [5] (see also [22]
for modern presentation) which after carrying out the thermodynamic limit,
goes into the Onsager formula [16,17].

\noindent The method we choose to proceed with, is an appropriate
generalization of the one used in [22] which consists there in reducing the
problem of finding of the partition function for the Ising model to
calculation of $\overline{{\rm Tr}}\left( P_1...P_s\right) $, where $%
\overline{{\rm Tr}}$ is the normalized trace while $P$'s are linear
combinations of $\gamma $ matrices - generators of usual Clifford algebra
naturally assigned to the lattice.

\noindent Then the observation that $\overline{{\rm Tr}}\left(
P_1...P_{2s}\right) $ is just a Pfaffian [3] of an antisymmetric matrix
formed with scalar products of $P$'s leads one to calculation of the
determinant from this very matrix.

\noindent The method proposed in [5,22] is purely algebraic, and though
probably not the shortest one, it lacks ambiguities of other methods and
uses well established, simple language of Clifford algebras. For other
aspects and connotations of such an approach see [9,10,11,12].

\smallskip\ 

\noindent Our paper is organized as follows:

\noindent In the second section we write down the $Z_n$-vector Potts model
in a form resembling (and generalizing!) the Ising model without external
field and then we represent the transfer matrix (an element of $%
C_{2p}^{\left( n\right) }$ !) as a sum of expressions proportional to ${\rm %
Tr}\left( \gamma _{i_1}...\gamma _{i_s}\right) $, where this time $\gamma $%
's are generalized $\gamma $ matrices.

\noindent Equivalently, the Hamiltonian for this Potts model can be looked
upon as an Action for the $Z_n$ chiral model on the square lattice.

\noindent In the third section we derive the formula for ${\rm Tr}\left(
\gamma _{i_1}...\gamma _{i_s}\right) $.

The last section is to supply an inevitable information on $C_k^{\left(
n\right) }$ algebras for the Reader's convenience as well as some
calculations avoided in the current text in order to make the presentation
more transparent.

\smallskip\ 

We would like to end this introduction by quotation from R. J. Baxter's book
(see [2] p. 454): ''{\it The only hope that occurs to me is just as Onsager
(1944) and Kaufman (1949) originally solved the zero-field Ising model by
using the algebra of spinor operators, so there may be similar algebraic
methods for solving the eight-vertex and Potts models}''.

\noindent Our suggestion then is that these very algebras are just
generalized Clifford algebras, and the presented paper is aimed to deliver
arguments in favor of that point of view.

\subsection*{II. The transfer matrix as a polynomial in $\gamma $'s}

In the following the transfer matrix $M$ for the $Z_n$-vector Potts model is
represented in a form of ''multi-sum'' of expressions proportional to ${\rm %
Tr}\left( \gamma _{i_1}...\gamma _{i_s}\right) $.

Let us assign to the set of states for this $Z_n$ vector Potts model on a $%
p\times q$ torus lattice ($p$ rows, $q$ columns), a set

\smallskip\ 

\noindent {\bf Def:}\quad \quad \quad \quad \quad $S=\left\{ \left(
s_{i,k}\right) =\left( p\times q\right) ;\,\,s_{i,k}\in Z_n\right\} \,,$%
\quad \quad $\diamondsuit $

\noindent where we have chosen a multiplicative realization for the cyclic
group $Z_n=\left( \omega ^l\right) _{l=0}^{n-1}$ and $s_{ik}$ denotes a
matrix element of the $\left( p\times q\right) $ matrix. Here, naturally, $%
\omega $ denotes the primitive root of unity.

\smallskip\ 

\noindent The total energy $E$ is then given by: 
$$
-\frac{E\left[ \left( s_{i,k}\right) \right] }{kT}=a\sum%
\limits_{i,k=1}^{p,q}\left(
s_{i,k}^{-1}s_{i,k+1}+s_{i,k+1}^{-1}s_{i,k}\right)
+b\sum\limits_{i,k=1}^{p,q}\left(
s_{i,k}^{-1}s_{i+1,k}+s_{i+1,k}^{-1}s_{i,k}\right) \eqno{(2.1)}\label{2.1} 
$$

\noindent while the partition function is defined to be: 
$$
Z=\sum\limits_{\left( s_{i,k}\right) \in S}\exp \left\{ -\frac{E\left[
\left( s_{i,k}\right) \right] }{kT}\right\} \eqno{(2.2)}\label{2.2} 
$$

One sees that the total energy of the system as represented by (2.1) is at
the same time - a generalization of its own $Z_2$ Ising case and - the
Action of the corresponding chiral model (for connection between lattice
gauge theory and spin systems see [6] and references therein).

\noindent The partition function could be written in terms of transfer
matrix and for that purpose we introduce the following notation:

\smallskip\ 

\noindent {\bf Notation:} 
$$
\begin{array}{l}
\vec s\cdot \vec s^{\prime }=\sum\limits_{i=1}^ps_is_i^{\prime }\,\,,\quad
\quad \quad \vec s_k=\left( 
\begin{array}{c}
s_{1,k} \\ 
s_{2,k} \\ 
\vdots \\ 
s_{p,k}
\end{array}
\right) \,,\quad \vec s_k^{*}=\left( 
\begin{array}{c}
s_{1,k}^{*} \\ 
s_{2,k}^{*} \\ 
\vdots \\ 
s_{p,k}^{*}
\end{array}
\right) \\ 
\left( s_{i,k}\right) =\left( \vec s_1,\vec s_2,...,\vec s_q\right)
\end{array}
\eqno{(2.3)}\label{2.3} 
$$

\noindent With the notation (2.3) adopted, the partition function $Z$ may be
now rewritten in a form 
$$
Z=\sum\limits_{\vec s_1,...,\vec s_q}\exp \left\{ a\sum\limits_{k=1}^q\left(
\vec s_k^{*}\cdot \vec s_{k+1}+\vec s_{k+1}^{*}\cdot \vec s_k\right)
+b\sum\limits_{k=1}^q\left( \vec s_k^{*}\cdot \sum\nolimits_1\vec s_k+\vec
s_k\cdot \sum\nolimits_1\vec s_k^{*}\right) \right\} ,\eqno{(2.4)}\label{2.4}
$$

\noindent after the natural periodicity conditions have been imposed.

\smallskip\ 

\noindent {\bf Periodicity conditions:} 
$$
\vec s_{q+1}=\vec s_1\,\,,\quad \left( \vec s_k\right) _1=\left( \vec
s_k\right) _{p+1}\,;\quad k=1,...,q\,\,,\eqno{(2.5)}\label{2.5} 
$$

\noindent where $\left( \vec x\right) _i$ denotes the $i$-th component of $%
\vec x$.

\noindent The matrix $\sum\nolimits_1$ is a $p\times p$ generalized Pauli
matrix with matrix elements $\delta _{i+1,j}$ , where $i,j\in Z_p^{^{\prime
}}=\left\{ 0,1,...,p-1\right\} $ and ''$+$'' is understood as the $%
Z_p^{^{\prime }}$ group action.

We introduce also the $\sigma _1$ generalized Pauli matrix, which is one of
the three $\sigma _1,\sigma _2,\sigma _3$ - playing the same role in
representing $C_{2p}^{\left( n\right) }$ generalized Clifford algebras as
the usual ones in representing the ordinary $C_{2p}^{\left( 2\right) }$
Clifford algebras via well known tensor products of $\sigma $ matrices [4].

It is now obvious that $Z$ may be represented as 
$$
Z={\rm Tr}\,M^q\,\,, \eqno{(2.6)}\label{2.6} 
$$

\noindent as we have 
$$
Z=\sum\limits_{\vec s_1,...,\vec s_q}M\left( \vec s_1,\vec s_2\right)
M\left( \vec s_2,\vec s_3\right) ...M\left( \vec s_q,\vec s_1\right) \,\,, 
$$

\noindent where matrix elements of the transfer matrix $M$ are given by: 
$$
M\left( \vec s,\vec s^{\prime }\right) =\exp \left\{ 2b\,{\rm Re}\left( \vec
s^{*}\cdot \sum\nolimits_1\vec s\right) \right\} \exp \left\{ 2a\,{\rm Re}%
\left( \vec s^{*}\cdot \vec s^{\prime }\right) \right\} \,\,.\eqno{(2.7)}%
\label{2.7} 
$$

\noindent It is convenient to consider the matrix $M$ as a product $M=BA$,
where the corresponding matrix elements are identified as 
$$
\begin{array}{l}
\,\,\,A\left( \vec s^{\prime \prime },\vec s^{\prime }\right) =\exp \left\{
2a\,{\rm Re}\left( \vec s^{\prime \prime *}\cdot \vec s^{\prime }\right)
\right\} \\ 
{\rm and} \\ 
\,\,\,B\left( \vec s,\vec s^{\prime \prime }\right) =\exp \left\{ 2b\,{\rm Re%
}\left( \vec s^{*}\cdot \sum\nolimits_1\vec s\right) \right\} \delta \left(
\vec s,\vec s^{\prime \prime }\right) \,\,.
\end{array}
\eqno{(2.8)}\label{2.8} 
$$

\noindent As all these $A,B,M$ matrices are multiindexed it is obvious that
they might be represented either as tensor products of $\left( n\times
n\right) $ matrices ($p$ times) or as $\left( n^p\times n^p\right) $
matrices.

\noindent It is not difficult then to see that 
$$
A=\otimes ^p\widehat{a}\quad \quad \quad {\rm i.e.} \eqno{(2.9)}\label{2.9} 
$$

\noindent $A$ is the $p$-th tensor power of the $\left( n\times n\right) $
matrix $\widehat{a}$, which has the form of a circulant matrix $W\left[
\sigma _1\right] $: 
$$
\widehat{a}=\left( \widehat{a}_{I,J}\right) =\left( \exp \left\{ 2a\,{\rm Re}%
\left( \omega ^{J-I}\right) \right\} \right) =\sum\limits_{l=0}^{n-1}\lambda
_l\sigma _1^l\equiv W\left[ \sigma _1\right] \,\,,\eqno{(2.10)}\label{2.10} 
$$

\noindent where $I,J\in Z_n^{^{\prime }}=\left\{ 0,1,2,...,n-1\right\} $ and 
$$
\lambda _l=\exp \left\{ 2a\,{\rm Re}\left( \omega ^l\right) \right\} \,\,.%
\eqno{(2.11)}\label{2.11} 
$$

In order to see that $A$ and $B$ matrices are just some elements of $%
C_{2p}^{\left( n\right) }$ we shall express them in terms of operators $X_k$
and $Z_k$\thinspace ; $k=1,2,...,p$ i.e. matrices typical for tensor product
representation of generalized Clifford algebras via generalized Pauli
matrices (see (A.3)).

\smallskip\ 

\noindent {\bf Def:} 
$$
\begin{array}{l}
X_k=I\otimes ...\otimes I\otimes \sigma _1\otimes I\otimes ...\otimes I\quad
\quad \quad \quad {\rm (}p{\rm -terms),} \\ 
Z_k=I\otimes ...\otimes I\otimes \sigma _3\otimes I\otimes ...\otimes I\quad
\quad \quad \quad {\rm (}p{\rm -terms),}\quad \diamondsuit
\end{array}
$$

\noindent where $\sigma _1$ and $\sigma _3$ are situated on the $k$-th site,
counting from the left hand side.

\smallskip\ 

\noindent The matrix $A$ may be therefore now rewritten as a product of $%
\left( n^p\times n^p\right) $ matrices 
$$
A=\prod\limits_{k=1}^pW\left[ X_k\right] \,\,,\quad \quad {\rm where}\quad
W\left[ X_k\right] =\sum\limits_{l=0}^{n-1}\lambda _lX_k^l\,\,\,.%
\eqno{(2.12)}\label{2.12} 
$$

\noindent Similarly, for the matrix $B$ we derive: 
$$
B=\exp \left\{ b\sum\limits_{k=1}^p\left(
Z_k^{-1}Z_{k+1}+Z_{k+1}^{-1}Z_k\right) \right\} \,\,\,, \eqno{(2.13)}\label
{2.13} 
$$

\noindent where $Z_{p+1}=Z_1$.

\smallskip\ 

\noindent The formula (2.13) follows from the simple observation that matrix
elements of \\$Z_k^{-1}Z_{k+1}+Z_{k+1}^{-1}Z_k$ (multiindexed by $%
\overrightarrow{s}$ and $\overrightarrow{s^{\prime \prime }}$) give exactly
ln of the corresponding term of (2.8) expression for $B$.

\noindent The $\delta $ function arises due to the fact that $\sigma
_3=\left( \delta _{I,J}\omega ^I\right) $ (see Appendix) and the
exponentiation of matrix elements is possible because $B$ is simply
proportional to unit matrix.

Once $A$ and $B$ have been represented as in (2.12) and (2.13) it is easy to
express them in terms of generalized $\gamma $ matrices. Introducing then
the tensor product representation (A.3) we get: 
$$
\begin{array}{l}
X_k=\omega ^{n-1}\gamma _k^{n-1}\overline{\gamma }_k \\ 
Z_k^{-1}Z_{k+1}=\overline{\gamma }_k^{n-1}\gamma _{k+1}\quad \quad {\rm %
for\,\,odd\,\,}n\,,
\end{array}
\eqno{(2.14)}\label{2.14} 
$$

\noindent and 
$$
\begin{array}{l}
X_k=\xi \omega ^{n-1}\gamma _k^{n-1}\overline{\gamma }_k \\ 
Z_k^{-1}Z_{k+1}=\xi \overline{\gamma }_k^{n-1}\gamma _{k+1}\quad \quad {\rm %
for\,\,even\,\,}n\,,
\end{array}
\eqno{(2.15)}\label{2.15} 
$$

\noindent where $k=1,2,...,p-1$ and $\xi ^2=\omega $ .

\noindent The corresponding expression on the boundaries - read: 
$$
Z_p^{-1}Z_1=U\overline{\gamma }_p^{n-1}\gamma _1\quad \quad {\rm %
for\,\,odd\,\,}n\,, \eqno{(2.16)}\label{2.16} 
$$

\noindent and 
$$
Z_p^{-1}Z_1=\xi ^{-1}U\overline{\gamma }_p^{n-1}\gamma _1\quad \quad {\rm %
for\,\,even\,\,}n\,, \eqno{(2.17)}\label{2.17} 
$$

\noindent where 
$$
\omega \cdot U=\otimes ^p\sigma _1\,.\eqno{(2.18)}\label{2.18} 
$$

\noindent For the proof of (2.14)-(2.17) use (A.6) and (A.7).

\smallskip\ 

\noindent From now on we shall proceed with formulas for $n$ being odd,
loosing nothing from generality of considerations while corresponding
formulas for the case of $n$ even are easily derivable from those for the
odd case. This in mind we get 
$$
A=\prod\limits_{k=1}^pW\left[ \omega ^{-1}\gamma _k^{n-1}\overline{\gamma }%
_k\right] \,\,,\eqno{(2.19)}\label{2.19} 
$$

$$
B=\exp \left\{ b\sum\limits_{k=1}^{p-1}\left( \overline{\gamma }%
_k^{n-1}\gamma _{k+1}+\gamma _{k+1}^{n-1}\overline{\gamma }_k\right)
\right\} \times \exp \left\{ bU\overline{\gamma }_p^{n-1}\gamma
_1+bU^{-1}\gamma _1^{n-1}\overline{\gamma }_p\right\} \,\,. \eqno{(2.20)}%
\label{2.20} 
$$

\noindent Our first goal is then achieved if one notes that 
$$
U=\prod\limits_{k=1}^p\gamma _k^{n-1}\overline{\gamma }_k \eqno{(2.21)}\label
{2.21} 
$$

\noindent i.e. the transfer matrix $M$ is now expressed in terms of
generalized $\gamma $ \thinspace matrices.

\noindent Before proceed, it is rather trivial and important to note that $%
U^n={\bf 1}$, $Z_k^n={\bf 1}$, $X_k^n={\bf 1}$, with obvious implication of
the same property for the $n$-th order polynomials in (2.19) and (2.20).

Our second and the main goal of this section is to represent the transfer
matrix in a form - reducing the ${\rm Tr}\,M^q$ problem to calculation of $%
{\rm Tr}\left( \gamma _{i_1}\gamma _{i_2}...\gamma _{i_s}\right) $ for some
collections of $\gamma $'s.

\noindent (Note that for $n=2$ the way to get the complete partition
function is shorter as there, it is enough to reduce the ${\rm Tr}\,M^q$
problem to calculation of ${\rm Tr}\left( P_1P_2...P_s\right) $ where $P$'s
are linear combinations of $\gamma $'s. Hence the number of necessary
summations is much, much smaller than in the case $n>2$, where it is rather
useless to try to represent $A$ and $B$ matrices in that convenient form).

\noindent For that to do we shall deal now with the matrix $B$, to reveal
its adequate for the purpose - structure.

\noindent Let us start with observations following from $U^n=1$ property of $%
U$. Define $V_k^{\pm }$ matrices (see (A.11)) to be 
$$
\begin{array}{l}
V_k^{+}=\frac 1n\sum\limits_{i=0}^{n-1}\omega ^{-ki}U^i \\ 
V_k^{-}=\frac 1n\sum\limits_{i=0}^{n-1}\omega ^{-ki}U^{-i}
\end{array}
\eqno{(2.22)}\label{2.22} 
$$

\noindent and also 
$$
\begin{array}{l}
\widetilde{B}_k^{+}=\exp \left\{ b\,\omega ^k\,\overline{\gamma }%
_p^{n-1}\gamma _1\right\} \\ 
\widetilde{B}_k^{-}=\exp \left\{ b\,\omega ^k\,\gamma _1^{n-1}\overline{%
\gamma }_p\right\} \,\,,\quad \quad {\rm where}\quad k=0,1,...,n-1\,.
\end{array}
\eqno{(2.23)}\label{2.23} 
$$

\noindent Then we have 
$$
\exp \left\{ b\,U\,\overline{\gamma }_p^{n-1}\gamma _1\right\}
=\sum\limits_{k=0}^{n-1}\widetilde{B}_k^{+}V_k^{+}\eqno{(2.24)}\label{2.24} 
$$

\noindent and 
$$
\exp \left\{ b\,U^{-1}\,\gamma _1^{n-1}\overline{\gamma }_p\right\}
=\sum\limits_{k=0}^{n-1}\widetilde{B}_k^{-}V_k^{-}\,\,\,.\eqno{(2.25)}\label
{2.25} 
$$

\noindent Equations (2.24) and (2.25) follow from examining $\exp \left\{
Ux\right\} $ via series expansions modulo $n$, as in (A.8)-(A.10).

The matrices $V_k^{\pm }$ have an important property (see (A.11)) 
$$
\left[ V_k^{\pm }\right] ^n=V_k^{\pm } \eqno{(2.26)}\label{2.26} 
$$

\noindent hence from (2.20), (2.24), (2.25) and the commutativity of matrix
arguments of $B$\\($U$, $\overline{\gamma }_k^{n-1}\gamma _{k+1}$...
etc.) it follows that: 
$$
B=\sum\limits_{l,k=0}^{n-1}\left( B_k^{+}V_k^{+}B_l^{-}V_l^{-}\right) %
\eqno{(2.27)}\label{2.27} 
$$

\noindent where 
$$
B_k^{+}=\exp \left\{ b\sum\limits_{\alpha =1}^{p-1}\overline{\gamma }_\alpha
^{n-1}\gamma _{\alpha +1}\right\} \widetilde{B}_k^{+}\,\,\,, \eqno{(2.28)}%
\label{2.28} 
$$

\noindent and 
$$
B_k^{-}=\exp \left\{ b\sum\limits_{\alpha ^{\prime }=1}^{p-1}\gamma _{\alpha
^{\prime }+1}^{n-1}\overline{\gamma }_{\alpha ^{\prime }}\right\} \widetilde{%
B}_k^{-}\,\,\,. \eqno{(2.29)}\label{2.29} 
$$

\noindent Expression (2.27) for $B$ becomes still simpler due to the
remarkable property of $V_k$'s: 
$$
V_kV_l=0\,\,\,;\quad \quad k\neq l \eqno{(2.30)}\label{2.30} 
$$

\noindent where, for the moment, $V_k=V_k^{\pm }$, (see Appendix for proof).

\noindent Hence 
$$
B=\sum\limits_{k=0}^{n-1}B_k^{+}B_k^{-}V_k^{+}V_k^{-}\,\,\,,\eqno{(2.31)}%
\label{2.31} 
$$

\noindent - all terms of the $k$-th summand - commuting.

As for the commuting of various matrices involved in representing of the
transfer matrix, note that 
$$
\left[ A,B\right] \neq 0\,\,,\quad \left[ U,A\right] =0\quad {\rm and}\quad
\left[ V_k^{+},A\right] =\left[ V_k^{-},A\right] =0\,\,. 
$$

\noindent All this is sufficient to write: 
$$
M^q=\sum\limits_{l=0}^{n-1}\left[ B_l^{+}B_l^{-}A\right] ^q\left(
V_l^{+}V_l^{-}\right) ^r\eqno{(2.32)}\label{2.32}
$$

\noindent where $1\leq r\leq n-1$ and $V^q=V^r$ as $V^n=V$.

It is not difficult to see how ''$r$'' arises. Namely ''$r$'' is the
residual of the quotient $\frac{\left( q-n\right) }{\left( n-1\right) }$.

This $n-1$ multiplicity in formula (2.32) makes the problem of thermodynamic
limit more interesting and involved.

\noindent If one assumes however that

\smallskip\ 

\noindent {\bf assumption:} 
$$
q=n+l\left( n-1\right) \,\,, \eqno{(2.33)}\label{2.33} 
$$

\noindent with $l$ an arbitrary integer, then 
$$
M^q=\sum\limits_{l=0}^{n-1}\left[ B_l^{+}B_l^{-}A\right]
^qV_l^{+}V_l^{-}\,\,\,. \eqno{(2.34)}\label{2.34} 
$$

\noindent Note that there is no multiplicity in (2.32) for the $Z_2$ case
(Ising model) and also for that case $V_k^{+}=V_k^{-}\equiv V_k$ , $k=0,1$.

The more: $V_kV_k=V_k$ and (2.34) reduces its form considerably.

\smallskip\ 

\noindent {\bf Comparison:}

In order to compare (2.34) with the known expression for $M^q$ in the case
of Ising model [22] one should note that, for $Z_2$ , (2.15) and (2.17)
differ only by $i$ and $-i$ from (2.14) and (2.16) correspondingly, hence we
have 
$$
M_q=\left( B_{-}A\right) ^qV_{+}+\left( B_{+}A\right) ^qV_{-}\,\,\,, %
\eqno{(2.35)}\label{2.35} 
$$

\noindent where 
$$
\begin{array}{l}
B_{-}=B_0^{+}B_0^{-}=\exp \left\{ 2bi\left( \sum\limits_{\alpha =1}^{p-1}%
\overline{\gamma }_\alpha \gamma _{\alpha +1}-\overline{\gamma }_p\gamma
_1\right) \right\} \,\,, \\ 
B_{+}=B_1^{+}B_1^{-}=\exp \left\{ 2bi\left( \sum\limits_{\alpha ^{\prime
}=1}^{p-1}\overline{\gamma }_{\alpha ^{\prime }}\gamma _{\alpha ^{\prime
}+1}-\overline{\gamma }_p\gamma _1\right) \right\}
\end{array}
$$

\noindent and 
$$
V_{+}=V_0=\frac 12\left( {\bf 1}+U\right) \,\,,\quad V_{-}=V_1=\frac
12\left( {\bf 1}-U\right) \,\,. \eqno{(2.35a)}\label{2.35a} 
$$

\noindent The formula (2.35) coincides then with the one known for Ising
model [22] apart from the obvious (see (2.1)) and insignificant scaling of
constants $a$ and $b$ by factor 2.

\noindent Using the formula (2.35), the notion of Pfaffian and its relation
to determinant - the author of [22] reobtained the complete partition
function leading to the famous Onsager formula.\quad $\diamondsuit $

\smallskip\ 

\noindent Having the same goal in mind we are going at first to examine the
expression (2.34) in order to see how (the polynomial in $\gamma $'s !) $M^q$
is represented as a multisum of summands proportional to ${\rm Tr}\left(
\gamma _{i_1}...\gamma _{i_s}\right) $. For that purpose we write: 
$$
\widetilde{B}_k^{+}=\exp \left( b\rho u_k^{+}\right) \,,\eqno{(2.36)}\label
{2.36} 
$$

\noindent and 
$$
\widetilde{B}_k^{-}=\exp \left( b\rho ^{-1}u_k^{-}\right) \,,\eqno{(2.37)}%
\label{2.37} 
$$

\noindent where 
$$
{\rm !}\quad \left( u_k^{+}\right) ^n=\left( u_k^{-}\right) ^n=1 
$$

\noindent i.e. 
$$
u_k^{+}\equiv \omega ^k\rho ^{-1}\overline{\gamma }_p^{n-1}\gamma _1\,\,\,, %
\eqno{(2.38)}\label{2.38} 
$$

$$
u_k^{-}\equiv \omega ^k\rho \gamma _1^{n-1}\overline{\gamma }_p\,\,\,, %
\eqno{(2.39)}\label{2.39} 
$$

\noindent while 
$$
\rho \equiv \rho \left( n\right) =\omega ^{\frac{n^2-1}2}\quad \quad \left(
k=0,1,...,n-1\right) .\eqno{(2.40)}\label{2.40} 
$$

\noindent Therefore both $\widetilde{B}_k^{+}$ and $\widetilde{B}_k^{-}$
become $n-1$ order polynomials in $u_k^{+}$ and $u_k^{-}$ correspondingly
(see (A.10)).

\noindent One also shows easily that $v_\alpha ^{+}$ and $v_\alpha ^{-}$
defined below

$$
\exp \left\{ b\,\overline{\gamma }_\alpha ^{n-1}\,\gamma _{\alpha
+1}\right\} =\exp \left\{ b\,\rho \,v_\alpha ^{+}\right\} \,\,\,,%
\eqno{(2.41)}\label{2.41} 
$$

$$
\exp \left\{ b\,\gamma _{\alpha +1}^{n-1}\,\overline{\gamma }_\alpha
\right\} =\exp \left\{ b\rho ^{-1}v_\alpha ^{-}\right\} \,\,\,,\eqno{(2.42)}%
\label{2.42} 
$$

\noindent do satisfy: 
$$
\left( v_\alpha ^{+}\right) ^n=\left( v_\alpha ^{-}\right) ^n=1\quad \quad
\alpha =1,...,p-1\,\,. \eqno{(2.43)}\label{2.43} 
$$

\noindent hence both expressions (2.41) and (2.42) become $n-1$ order
polynomials in matrices $v_\alpha ^{+}$ and $v_\alpha ^{-}$ correspondingly.

\smallskip\ 

\noindent {\bf n=2:}

\noindent For $n=2$ the further job is extremely facilitated due to the fact
that $A$ becomes then of the form 
$$
A=\prod\limits_{k=1}^pP_kQ_k \eqno{(2.44)}\label{2.44} 
$$

\noindent where $P$'s and $Q$'s are some known linear combinations of
ordinary $\gamma$ matrices and similar holds for $B_{+}$, $B_{-}$ matrices
from (2.35).

\noindent The matrices $V_{+}$, $V_{-}$ from (2.35) have also simple form
and thus ${\rm Tr}\,M^q$ disentangles for $n=2$ to be the sum only four
summands of the Pfaffian type i.e. ${\rm Tr}\left( P_1P_2...P_s\right) $.

\smallskip\ 

\noindent {\bf n}${\bf >}${\bf 2:}

\noindent Unfortunately this disentanglement is no more possible for $n>2$
as the polynomial $W$ (see (2.12)) no more is representable uniquely as the
product of linear combinations of $\gamma $'s and neither is $B$.

\noindent Hence the multi-sum becomes more complicated.

\noindent Nevertheless it is obvious that ${\rm Tr}\,M^q$ problem reduces to
calculation of ${\rm Tr}\left( \gamma _{i_1}\gamma _{i_2}...\gamma
_{i_s}\right) $ for some collections of $\gamma $'s.

\smallskip\ 

\noindent {\bf n=2:}

\noindent In the case of Ising model, the four arising Pfaffians contribute
to the partition function to give [22]: 
$$
\begin{array}{l}
Z=2^{pq-1}\left\{ \prod\limits_{k,l=1}^{p,q}\left[ {\rm ch}2a^{\prime }{\rm %
ch}2b^{\prime }-{\rm sh}2a^{\prime }\cos \frac \pi q\left( 2l+1\right) -{\rm %
sh}2b^{\prime }\cos \frac \pi p\left( 2k+1\right) \right] ^{\frac
12}+\right.  \\ 
\quad \quad +\prod\limits_{k,l=1}^{p,q}\left[ {\rm ch}2a^{\prime }{\rm ch}%
2b^{\prime }-{\rm sh}2a^{\prime }\cos \frac \pi q\left( 2l+1\right) -{\rm sh}%
2b^{\prime }\cos \frac{2\pi k}p\right] ^{\frac 12}+ \\ 
\quad \quad +\prod\limits_{k,l=1}^{p,q}\left[ {\rm ch}2a^{\prime }{\rm ch}%
2b^{\prime }-{\rm sh}2a^{\prime }\cos \frac{2\pi l}q-{\rm sh}2b^{\prime
}\cos \frac \pi p\left( 2k+1\right) \right] ^{\frac 12}+ \\ 
\quad \quad \left. -\sigma \prod\limits_{k,l=1}^{p,q}\left[ {\rm ch}%
2a^{\prime }{\rm ch}2b^{\prime }-{\rm sh}2a^{\prime }\cos \frac{2\pi l}q-%
{\rm sh}2b^{\prime }\cos \frac{2\pi k}p\right] ^{\frac 12}\right\} 
\end{array}
\eqno{(2.45)}\label{2.45}
$$

\noindent where $\sigma $ denotes the sign of $T-T_c$ and $a^{\prime }=2a$, $%
b^{\prime }=2b$. Both the square root and the $\sigma $-sign have appeared
here because of the use of $Pf^2=\det $ relation.

\smallskip\ 

\noindent {\bf n}${\bf >}${\bf 2:}

\noindent Again, for $n>2$, as we shall see in the following section,
although the generalization of the Pfaffian is possible to the case of
arbitrary $n$, its relation to any valuable generalization of determinant
does not to be valid as the arising signum like function no more is an
epimorphism of $S_k$ onto $Z_n$, except for $n=2$ of course. ($S_k$ - the
symmetric group of $k$-elemental permutations). However, one may write $M^q$
as the polynomial in $\gamma $'s and then use the general formula from the
following section.

\noindent This representation of $M^q$ in terms of $\gamma _{i_1}...\gamma
_{i_s}$ products is given in the Appendix.

\noindent As a result we have the following structure of the complete
partition function for the $Z_n$ vector Potts models: 
$$
M^q=\frac 1{n^2}\sum\limits_{j_1,j_2=0}^{n-1}\sum\limits_{\vec \Pi \in
\Gamma _m}G\left( \vec \Pi \right) \widehat{\Omega }\left( \vec \Pi
;j_1,j_2\right) \,\,,\eqno{(2.46)}\label{2.46} 
$$

\noindent where $G\left( ...\right) $ are {\bf known functions} of
parameters $a$ and $b$, (see the Appendix: (A.15)) and 
$$
\overline{{\rm Tr}}{\rm \,}\widehat{\Omega }=\omega ^i\quad or\quad {\rm Tr}%
=0\,\,, 
$$
\noindent where $i=i\left( \vec \Pi ,j_1,j_2\right) \in Z_n^{\prime
}=\left\{ 0,1,...,n-1\right\} $.

The dependence of $i$ on its indices is easy to be derived using the most
general, appropriate formula for $\overline{{\rm Tr}}\left( \gamma
_{i_1}\gamma _{i_2}...\gamma _{i_s}\right) $ supplied by the next section.

\medskip\ 

\subsection*{III. Trace formula for any element of $C_{2p}^{\left( n\right)
} $}

In this section the explicit formula for trace of any element of $%
C_{2p}^{\left( n\right) }$ algebra is delivered.

\noindent The very formula is crucial for getting the complete partition
function for Potts models and hence (see [2] p. 454) for solving several
major problems of statistical physics being unsolved till now since many
years.

The problem of explicit trace formula for $M^q\in C_{2p}^{\left( n\right) }$%
, is decisive in calculation of $Z$ function for those models on the lattice
in which the transfer matrix is an element of $C_{2p}^{\left( n\right) }$.

\noindent We proceed now to derivation of the very formula.

\smallskip\ 

\noindent Note! By definition, in this section Tr map is normalized i.e. $%
{\rm Tr}\,I=1$. The derivation has the form of a sequence of lemmas.

\smallskip\ 

\noindent {\bf Lemma 1.}

Let 
$k\neq n\,\,\,\,mod\,n\,,\,\,\,k\in N\,$; \thinspace \thinspace 
then ${\rm Tr}\left( \gamma _{i_1}...\gamma _{i_k}\right) =0$
.\quad \quad $\diamondsuit $

\smallskip\ 

\noindent Proof: The same as for usual Clifford algebras. Use the matrix $U$
defined by (2.18).\quad \quad $\diamondsuit $

\smallskip\ 

\noindent {\bf Lemma 2.}

${\rm Tr}\left( \gamma _{i_1}\gamma _{i_2}...\gamma _{i_k}\right) \neq 0$ 
{\bf iff} there exists permutation $\delta \in S_{kn}$ , such that

\noindent
$i_{\sigma \left( 1\right) }=i_{\sigma \left( 2\right) }=...=i_{\sigma
\left( n\right) }$,\quad $i_{\sigma \left( n+1\right) }=...=i_{\sigma \left(
2n\right) }$,\quad ...\quad ,\quad $i_{\sigma \left( kn-n+1\right)
}=...=i_{\sigma \left( kn\right) }$.\quad $\diamondsuit $

\smallskip\ 

\noindent Proof: The proof follows from observation that due to (A.1) if no $%
n$-tuple of the same $\gamma $'s exists then ${\rm Tr}\left( ...\right) =0$.
Other steps of the proof are reduced to this first one.\quad \quad $%
\diamondsuit $

\smallskip\ 

\noindent It is therefore trivial to note, but important to realize, that:

\smallskip\ 

\noindent {\bf Lemma 3.}

${\rm Tr}\left( \gamma _{i_1}...\gamma _{i_k}\right) =0$ or $l\in Z_n$ - the
multiplicative cyclic group of $n$-th roots of unity.\quad \quad $%
\diamondsuit $

\smallskip\ 

\noindent In Lemma 3 $k$ is again an arbitrary integer while in all
preceding lemmas, and in the following, $i_1,i_2,...,i_k$ run from 1 to
number of generators of the given algebra. This number was chosen to be
even, however note [15] that the ''odd case'' problem is reduced to this
very one due to the properties of generalized Clifford algebra
representations.

The major problem now is to determine this value ''$0$ or $l$'' for
arbitrary set of indices $i_1,i_2,...,i_k$.

\noindent In order to do that define a signum like function $K$
(unfortunately it is an epimorphism only for $n=2$) - as follows:

\smallskip\ 

\noindent {\bf Def:} 
$$
K:S_p\rightarrow Z_n\,\,;\quad \quad \quad \Theta _{\sigma \left( 1\right)
}\Theta _{\sigma \left( 2\right) }...\Theta _{\sigma \left( p\right)
}=K\left( \sigma \right) \Theta _1\Theta _2...\Theta _p\quad , 
$$
where $\Theta $'s satisfy (A.1) except for the condition $\gamma _i^n=1$,
which is now replaced by $\Theta _i^2=1$.\quad \quad $\diamondsuit $

\smallskip\ 

\noindent This definition being adapted, it is now not very difficult to
prove:

\smallskip\ 

\noindent {\bf Lemma 4.}

${\rm Tr}\left( \gamma _{i_1}...\gamma _{i_{pn}}\right) =K\left( \Sigma
\right) K\left( \sigma \right) $ ,\quad \quad for

\noindent a)\quad $i_{\sigma \left( 1\right) }=...=i_{\sigma \left( n\right)
}\,\,$,\quad ...\quad ,\quad $i_{\sigma \left( pn-n+1\right) }=...=i_{\sigma
\left( pn\right) }$\quad and

\noindent b)\quad $i_{\widetilde{\sigma }\left( n\right) }\,<i_{\widetilde{%
\sigma }\left( 2n\right) }\,<...<i_{\widetilde{\sigma }\left( pn\right) }\,$,

\noindent where $\widetilde{\sigma }\equiv \Sigma \circ \sigma $ , while $%
\Sigma $ is a permutation of the elements $\left\{ n,2n,...,pn\right\} $.
(The group of $\Sigma $'s is naturally identified with an appropriate
subgroup of $S_{pn}$).\quad \quad $\diamondsuit $

\smallskip 

\noindent Proof: The proof relies on observation that these are only {\bf %
different} $n$-tuples which are ''rigidly'' shifted ones trough the others,
i.e. there is no permutation within any given $n$-tuple.\quad \quad $%
\diamondsuit $

\smallskip 

\noindent The generalization of Lemma 4 to the arbitrary case of some of the 
$n$-tuples being equal - is straightforward. (The necessary change of
conditions a) and b) is obvious).

This in mind and from other lemmas we finally get:

\medskip\ 

\noindent {\bf Theorem:}

${\rm Tr}\left( \gamma _{i_1}...\gamma _{i_{pn}}\right) =\underset{\sigma
\in S_{pn}}{\sum^{\prime }}\sum\limits_{\vec p}\sum\limits_{\Sigma \in
S_{\vec p}}K\left( \Sigma \right) K\left( \sigma \right) \times \delta
\left( i_{\widetilde{\sigma }\left( 1\right) },...,i_{\widetilde{\sigma }%
\left( p_1n\right) }\right) \times $

$\quad \quad \times \delta \left( i_{\widetilde{\sigma }\left( p_1n+1\right)
},...,i_{\widetilde{\sigma }\left( \left[ p_1+p_2\right] n\right) }\right)
\times ...\times \delta \left( i_{\widetilde{\sigma }\left( pn-p_ln+1\right)
},...,i_{\widetilde{\sigma }\left( pn\right) }\right) $ ,

\noindent with the notation to follow.\quad \quad $\diamondsuit $

\smallskip 

\noindent {\bf Notation:}

$\vec p=\left( p_1,p_2,...,p_l\right) $, $p_i\geq 1$, $\sum%
\limits_{i=1}^lp_i=p$, $\widetilde{\sigma }=\Sigma \circ \sigma $ , and $%
S_{\vec p}$ is a subgroup of $S_{pn}$ isomorphic (for example!) to the group
of all block matrices obtained via permutations of ''block columns'' of the
matrix

$\left( 
\begin{array}{cccc}
I_{p_1\!^n} &  &  &  \\ 
& I_{p_2\!^n} &  &  \\ 
&  & \ddots &  \\ 
&  &  & I_{p_l\!^n}
\end{array}
\right) $ , \quad where $I_k$ is the $\left( k\times k\right) $ unit matrix.

\noindent $\delta $ - here denotes the multi-indexed Kronecker delta i.e. it
assigns zero unless all its arguments are equal and in this very case $%
\delta \left( ...\right) =1$. The sum $\Sigma ^{\prime }$ is meant to take
into account {\bf only} those permutations that do satisfy the conditions:

a)\quad $\sigma \left( 1\right) <\sigma \left( 2\right) <...<\sigma \left(
p_1n\right) $,\quad ...\quad , $\sigma \left( pn-p_ln+1\right) <...<\sigma
\left( pn\right) $,

and

b)\quad $\sigma \left( 1\right) <\sigma \left( p_1n+1\right) <...<\sigma
\left( pn-p_ln+1\right) $.

\smallskip 

\noindent {\bf Comments:}

\noindent 1) For the case of $n=2$ the theorem gives us the Pfaffian of the
product $\gamma _{i_1},...,\gamma _{i_{p^2}}$, as in the case, (and only!
for $n=2$) $K\left( \Sigma \right) =1$ and we are left, as a result with
only $\Sigma ^{\prime }$ sum, while Kronecker deltas become functions of the
same number of indices $i_j$.

\noindent 2) The theorem solves our problem of ${\rm Tr}\,M^q$, as any
element of generalized Clifford algebra is a polynomial in $\gamma $'s
satisfying (A.1), including $M^q\in C_{2p}^{\left( n\right) }$.

\subsection*{IV. Final comments for the Part 1 of the presentation}

We have carried out our twentieth century investigation for the $Z_n$ vector Potts model known
also under the name of planar Potts model.

\noindent The similar investigation of the other Potts models, i.e. standard
Potts models with two-site interaction [23] and multisite interactions as
well, is being now carried out.

However, it is to be noted here that the model considered in [14] possesses
transfer matrix $M=BA$, where matrix $A$ is a particular case of the one
defined by (2.19) while $B$, though also expressed by $Z_k$ matrices defined
in section II, has a different polynomial (in these operators) in the
expotential.

As for the multisite interactions, the algebras to be used are the universal
generalized Clifford algebras, introduced in [7]. Needless to say that these
are standard Potts models which are of more interest because of their
relation to a number of outstanding problems in lattice statistics [2,23].

\noindent To this end let us express our suggestion that the models of lattice
statistics could be adapted (thanks to specific interpretation) to the
domain of urban economics which, using the notion of entropy and information
introduces as a matter of fact a kind of thermodynamics [13].

\subsection*{V. An outline of the second planned part  content }

The content of the preceding chapters-except for the  algorithm for calculation of the complete partition function
as presented above - was already published in twenty first century [25],[26].
This especially  concerns the Clifford algebra  technique omnipresent here and planned to play
similar leading role in the next part of our review.
The content of [26] from 2001 indicates the idea and a way how to use generalized Clifford algebra for chiral Potts models on the plane.
The incessantly growing area of applications of Clifford algebras and naturalness of their use in formulating problems for direct calculation entitles one to call them Clifford numbers. The generalized ``universal'' Clifford numbers are here introduced via {k}-ubic form $Q_k$ replacing quadratic one in familiar construction of an appropriate ideal of tensor algebra. One of the epimorphic images of universal algebras $k-C_n \equiv T(V)/I(Q_k)$ is the algebra $Cl_n(k)$ with {n} generators and these are the algebras to be used here. Because generalized Clifford algebras $Cl_n(k)$ possess inherent $Z_k$ grading - this property makes them an efficient apparatus to deal with spin lattice systems. This efficiency is illustrated in [26] by derivation of two major observations. Firstly, the  partition functions for vector and planar Potts models and other model with $Z_n$ invariant Hamiltonian are polynomials in generalized hyperbolic functions of the {n}-th order. Secondly, the problem of algorithmic calculation of the partition function for any vector Potts model as treated here is reduced to the calculation of traces of products of the generators of the generalized Clifford algebra. Finally the expression for such traces for arbitrary collection of generator matrices is derived  in [26].

\noindent Since the same  2001 year, due to the authors of [27] we know the form of the $k$-state Potts model partition function (equivalent to the Tutte polynomial) for a lattice strip of fixed width and arbitrary length.

\noindent From  2005 year Alan D. Sokal 54 pages review  [28] aimed for mathematicians too one may learn that "`
the multivariate Tutte polynomial (known to physicists as the Potts-model partition function) can be defined on an arbitrary finite graph G, or more generally on an arbitrary matroid M, and encodes much important combinatorial information about the graph"'.

Alan D. Sokal discusses there "`some questions concerning the complex zeros of the multivariate Tutte polynomial, along with their physical interpretations in statistical mechanics (in connection with the Yang--Lee approach to phase transitions) and electrical circuit theory."'  Quite numerous open problems are also posed in [28]. For many references see  both [27] and [28].

\noindent Coming back for a while to  generalized Clifford algebra we mark their being just only mentioned  in the abstract of Baxter paper [29]
in which: "`The partition function of the N-state superintegrable chiral Potts model is obtained exactly and explicitly (if not completely rigorously) for a finite lattice with particular boundary conditions"'.  "` The associated Hamiltonian has a very simple form, suggesting that may be a more direct algebraic method (perhaps a generalized Clifford algebra) for obtaining its eigenvalues."'
Baxter - including his famous book [1982 Exactly Solved Models in Statistical Mechanics] comes back several time to the idea of Clifford or Clifford-like algebras' potential importance for the still not solved problem of complete partitions function obtained in a way Ising model was splved with help of Clifford algebras properties being  used in a natural and elegant way. Manageable in a understandable way.
Let us quote after Baxter from [30]:"' "This is rather intriguing - we are in much the same position with the chiral Potts model
as we were in 1951, when Professor Yang entered the field of statistical mechanics by calculating
MO for the Ising model. So on the occasion of his 70th birthday we are able to present him not
only with this meeting in honor of his great contributions to theoretical physics, but also with
an outstanding problem worthy of his mettle. Plus a change, c'est la $m\hat{e}me$ chose."'Here and there 

P. P. Martin's book  [31] and his many subsequent papers in twentieth as well in twenty first century [32] are inevitable source of ideas and inspiration.

\noindent And here now comes the 2005 year. 2005   was declared by some authors to be the major breakthrough for the chiral Potts model. This concerns 

an outstanding problem of the order parameters. See Baxter again  [33] and [34]. There in [33] and then in [34] Baxter deals again with 
the problem of the  order parameter in the chiral Potts model. He recalls that an elegant conjecture for this was made in 1983 and that it has since been successfully tested against series expansions, but as far as the author of [34] is aware there is as yet no proof of the conjecture.

2005. Again Professor  Baxter. Here is an abstract  of his Annual Conference 2005 of the  Australian Mathematical Society  public talk             .:
entitled \textbf{Lattice models in statistical mechanics: the chiral Potts model}:"`There are a few lattice models of interacting systems that can be solved exactly, in the sense that one can calculate the free energy in the thermodynamic limit of a
large system. The interesting ones are mostly two-dimensional, such as the Ising model and
the six and eight-vertex models. A comparitively recent addition to the list is the chiral Potts
model. This is more difficult mathematically than its predecessors. While its free energy
was calculated in 1988, until now there has only been a conjecture (a very elegant one) for
the order parameters, i.e. the spontaneous magnetizations. This conjecture has now been
verified, and in this talk I shall discuss the difficulties encountered and the method used. "'
The solution from 1988 he is referring to apparently refers to papers such as [35], [36], [37], [38] and others later, see [39].
By the way? - the authors of [39] implicitly indicate Generalized Clifford Algebras - in a footnote (3) referring there to
Morris papers from   1967 and 1968 [quoted in all Kwasniewski papers on subject].
The 69 pages, 30 figures reviw [39] was written in honor of Onsager's ninetieth birthday,also in order to present  "`some exact results  in the chiral Potts models and to translate these results into language more transparent to physicists"'.

This is more or less what the second part is planned to be about. By no means it sholud include review of
numerous cotributions of Professor  F.Y. Wu including not only Potts models [see:\\  
http://www.physics.neu.edu/Department/Vtwo/faculty/wu.../wupubupdated81803.htm \\
but such fascinating papers as [40] of specifically personal interest of the authors [http://ii.uwb.edu.pl/akk/publ1.htm].
Papers published in Advances in Applied Clifford Algebras such as [42]  and the papers by the authors are to be included in the second part of
this review too.

\subsection*{Appendix}

\noindent {\bf 1.}\quad $C_{2p}^{\left( n\right) }$ generalized Clifford
algebra is defined [15] to be generated by $\gamma _1,...,\gamma _{2p}$
matrices satisfying: 
$$
\gamma _i\gamma _j=\omega \gamma _j\gamma _i\,,\quad i<j\,,\quad \gamma _i^n=%
{\bf 1}\,,\quad i,j=1,2,...,2p\,. \eqno{(A.1)}\label{A.1} 
$$

\noindent The very algebra has - up to equivalence - only one irreducible
and faithful representation, and its generators can be represented as tensor
products of generalized Pauli matrices: 
$$
\sigma _1=\left( \delta _{i+1,j}\right) \,,\quad \sigma _2=\left( \omega
^i\delta _{i+1,j}\right) \,,\quad \sigma _3=\left( \omega ^i\delta
_{i,j}\right) \,, \eqno{(A.2)}\label{A.2} 
$$

\noindent where $i,j\in Z_n^{\prime }=\left\{ 0,1,...,n-1\right\} $ - the
additive cyclic group.

\noindent One easily checks, that $\left\{ \sigma _i\right\} _1^3$ do
satisfy (A.1) for $n$ being {\bf odd}.

\noindent Let $I$ denotes since now the unit $\left( n\times n\right) $
matrix and let 
$$
\begin{array}{l}
\gamma _1=\sigma _3\otimes I\otimes I\otimes ...\otimes I\otimes I\,\,, \\ 
\gamma _2=\sigma _1\otimes \sigma _3\otimes I\otimes ...\otimes I\otimes
I\,\,, \\ 
\vdots \\ 
\gamma _p=\sigma _1\otimes \sigma _1\otimes \sigma _1\otimes ...\otimes
\sigma _1\otimes \sigma _3\,\,, \\ 
\bar \gamma _1=\sigma _2\otimes I\otimes I\otimes ...\otimes I\otimes I\,\,,
\\ 
\bar \gamma _2=\sigma _1\otimes \sigma _2\otimes I\otimes ...\otimes
I\otimes I\,\,, \\ 
\vdots \\ 
\bar \gamma _p=\sigma _1\otimes \sigma _1\otimes \sigma _1\otimes ...\otimes
\sigma _1\otimes \sigma _2\,\,,
\end{array}
\eqno{(A.3)}\label{A.3} 
$$

\noindent then $\left\{ \gamma _i,\bar \gamma _j\,,\,\,i,j=1,...,p\right\} $
do satisfy (A.1) with $\omega $ replaced by $\omega ^{-1}$, hence (A.3) are
generators of the algebra isomorphic to $C_{2p}^{\left( n\right) }$
(isomorphism is given by $\sigma _1\leftrightarrow \sigma _3$ in (A.3))
(This very (A.3) representation was chosen for technical reason - we get,
for example, in calculations of section II, the matrix $U$ without
coefficients etc.).

\noindent It is also to be noted that for $n$ being {\bf odd} 
$$
\sigma _3=\sigma _1^{n-1}\sigma _2\,\,. \eqno{(A.4)}\label{A.4} 
$$

The case of $n$ being {\bf even} leads to similar representation with $%
\sigma _1$ unchanged but $\sigma _2$ and $\sigma _3$ now equal to: 
$$
\sigma _2=\left( \xi ^i\delta _{i+1,j}\right) \,,\quad \sigma _3=\xi \sigma
_1^{n-1}\sigma _2\,\,, \eqno{(A.5)}\label{A.5} 
$$

\noindent where $\xi $ is a {\bf primitive} $2n$-th root of unity such that 
$$
\xi ^2=\omega \,\,. 
$$

\noindent (A.3) then with these appropriate for case $n=2\nu $ generalized
Pauli matrices, provides us with the same type representation of $%
C_{2p}^{\left( n\right) }$ as the one for the case $n=2\nu +1$.

\noindent One then easily proves that 
$$
\begin{array}{l}
\sigma _3^{n-1}\sigma _2=\omega \sigma _1\,\,, \\ 
\sigma _2^{n-1}\sigma _1=\sigma _3^{-1}\quad \quad {\rm for}\quad n=2\nu +1
\end{array}
\eqno{(A.6)}\label{A.6} 
$$

\noindent and 
$$
\begin{array}{l}
\sigma _3^{n-1}\sigma _2=\xi ^{-1}\omega \sigma _1\,\,, \\ 
\sigma _2^{n-1}\sigma _1=\xi ^{-1}\sigma _3^{-1}\quad \quad {\rm for}\quad
n=2\nu \,\,.
\end{array}
\eqno{(A.7)}\label{A.7} 
$$

\medskip\ 

\noindent {\bf 2.}\quad In this part of the Appendix we derive one useful
formula, necessary for section II.

Let $x$ be any element of an associative, finite dimensional algebra with
unity ${\bf 1}$.

\noindent Then 
$$
\begin{array}{l}
\exp \left\{ x\right\} =\sum\limits_{i=0}^{n-1}f_i\left( x\right) \,\,,\quad
\quad {\rm where} \\ 
f_i\left( x\right) =\sum\limits_{k=0}^\infty \frac{x^{nk+i}}{\left(
nk+i\right) !}\,\,,\quad i=0,...,n-1\,.
\end{array}
\eqno{(A.8)}\label{A.8} 
$$

\noindent We now express these $f_i$'s in terms expotentials.

\noindent For that to do it is sufficient to note that 
$$
f_i\left( \omega x\right) =\omega ^if_i\left( x\right) \,\,,\quad \quad
i=0,1,...,n-1\,.\eqno{(A.9)}\label{A.9} 
$$

\noindent The (A.9) reveals the $Z_n$ symmetry properties of these
generalized ''cosh'' functions and we get from this set of relations 
$$
f_i\left( x\right) =\frac 1n\sum\limits_{k=0}^{n-1}\omega ^{-ki}\exp \left\{
\omega ^kx\right\} \,\,,\quad i=0,...,n-1\,. \eqno{(A.10)}\label{A.10} 
$$

\medskip\ 

\noindent {\bf 3.}\quad For considerations of the section II we need the
following

\smallskip\ 

\noindent {\bf Lemma}

Let $U$ be as $x$ above and in addition let $U^n={\bf 1}$. Then $V$ defined
as follows 
$$
\begin{array}{l}
V=\frac 1n\sum\limits_{i=0}^{n-1}U^i\,,\quad {\rm has\,\,the\,\,property:}
\\ 
V^n=V\,\,.\quad \quad \quad \quad \quad \quad \quad \quad \quad \diamondsuit
\end{array}
\eqno{(A.11)}\label{A.11} 
$$

\smallskip\ 

\noindent Proof: For the proof, just note that for some $a_i$'s 
$$
V^n=\sum\limits_{i=0}^{n-1}a_iU^i\,, 
$$

\noindent and both sides of this identity equation must be symmetric in $U^i$
monomials. One concludes therefore that $a_i=a_j$ , $i,j=0,1,...,n-1$, hence
- counting the number of all arising summands - one arrives at the
conclusion of the Lemma.

\medskip\ 

\noindent {\bf 4.}\quad Here, the proof of the (2.30) formula follows.

\noindent Let $k=l+r$, $0<r\leq n-1$. Then 
$$
V_k^{+}V_l^{+}=\frac
1{n^2}\sum\limits_{i_1=0}^{n-1}\sum\limits_{i_2=0}^{n-1}\omega
^{-ri_1}\left( \omega ^{-1}U\right) ^{i_1+i_2}. 
$$

\noindent Introduce now new summation indices: $i_1$ and $i=i_1+i_2$.

\noindent Then we have 
$$
V_k^{+}V_l^{+}=\frac 1{n^2}\sum\limits_{i_1=0}^{n-1}\omega
^{-ri_1}\sum\limits_{i=0}^{n-1}\left( \omega ^{-1}U\right) ^i=0\,, 
$$

\noindent because the summation over $i_1$ gives zero.

\noindent The proof for other $\left( +,-\right) $, $\left( -,+\right) $, $%
\left( -,-\right) $ cases is the same.

\medskip\ 

\noindent {\bf 5.}\quad In this part of the Appendix we derive the multisum
structure of the complete partition function $Z$, for the $Z_n$-vector Potts
model with any $n\geq 2$.

\noindent Since now on we shall use the following abbreviations: 
$$
\Gamma _k\equiv Z_n^{\prime }\otimes Z_n^{\prime }\otimes ...\otimes
Z_n^{\prime }\,\,,\quad {\rm (}k{\rm \,\,summands)}\,, 
$$

\noindent where $Z_n^{\prime }=\left\{ 0,1,...,n-1\right\} $, and 
$$
\omega ^{\frac{n^2-1}2}\equiv \rho \left( n\right) \equiv \rho \,\,. 
$$

\noindent Recall also $u_k^{+}$ , $u_k^{-}$ ; $k\in Z_n^{\prime }$ and $%
v_\alpha ^{+}$, $v_\alpha ^{-}$ ; $\alpha =1,...,p-1$, defined in section
II. Then we have, according to (A.10): 
$$
\begin{array}{l}
B_k^{+}=\sum\limits_{l_1=0}^{n-1}f_{l_1}\left( b\rho \right) \left(
u_k^{+}\right) ^{l_1}\prod\limits_{\alpha
=1}^{p-1}\sum\limits_{s=0}^{n-1}f_s\left( b\rho \right) \left( v_\alpha
^{+}\right) ^s\,\,, \\ 
B_k^{-}=\sum\limits_{l_2=0}^{n-1}f_{l_2}\left( b\rho ^{-1}\right) \left(
u_k^{-}\right) ^{l_2}\prod\limits_{\alpha ^{\prime
}=1}^{p-1}\sum\limits_{t=0}^{n-1}f_t\left( b\rho ^{-1}\right) \left(
v_\alpha ^{-}\right) ^t\,\,.
\end{array}
\eqno{(A.12)}\label{A.12} 
$$

In order to manage with abundance of indices we introduce a further,
deliberate notation.

\smallskip\ 

\noindent {\bf Notation:}

\quad \quad \quad \quad \quad $\vec L\in \Gamma _2\,\,\,,\quad \quad \quad
\quad \quad \vec T,\vec S\in \Gamma _{p-1}$\quad \quad \quad i.e.

\noindent $\vec L\equiv \left( l_1,l_2\right) $ where $l_1,l_2\in
Z_n^{\prime }$ etc.

\smallskip\ 

\noindent We also define: 
$$
\begin{array}{l}
F\left( b;k,\vec L,\vec S,\vec T\right) \equiv \omega ^{k\left(
l_1+l_2\right) }\rho ^{l_2-l_1}f_{l_1}\left( b\rho \right) f_{l_2}\left(
b\rho ^{-1}\right) \prod\limits_{i=1}^{p-1}f_{s_i}\left( b\rho \right) \rho
^{-s_i}\times \prod\limits_{j=1}^{p-1}f_{t_j}\left( b\rho ^{-1}\right) \rho
^{t_j}, \\ 
\widehat{\Xi }\left( \vec L,\vec S,\vec T\right) \equiv \left( \bar \gamma
_p^{n-1}\gamma _1\right) ^{l_1}\left( \gamma _1^{n-1}\bar \gamma _p\right)
^{l_2}\prod\limits_{i=1}^{p-1}\left( \bar \gamma _i^{n-1}\gamma
_{i+1}\right) ^{s_i}\prod\limits_{j=1}^{p-1}\left( \gamma _{j+1}^{n-1}\bar
\gamma _j\right) ^{t_j}
\end{array}
$$

\noindent and

\quad \quad $B_k\equiv B_k^{+}B_k^{-}$\thinspace \thinspace .\quad \quad $%
\quad \diamondsuit $

\smallskip\ 

\noindent Hence we may write 
$$
B_k=\sum_{\vec L,\vec S,\vec T}F\left( b;k,\vec L,\vec S,\vec T\right) 
\widehat{\Xi }\left( \vec L,\vec S,\vec T\right) \eqno{(A.13)}\label{A.13} 
$$

\noindent where 
$$
\vec L\in \Gamma _2\,\,,\quad \vec S,\vec T\in \Gamma _{p-1}\,\,. 
$$

\noindent It is important to recall now that: 
$$
M^q=\sum\limits_{k=0}^{n-1}\left[ B_k,A\right] ^qV_k^{+}V_k^{-}\quad \quad 
{\rm and}\quad \left[ B_k,A\right] \neq 0\,\,. 
$$

\noindent We introduce therefore again an appropriate notation.

\smallskip\ 

\noindent {\bf Notation:} 
$$
\begin{array}{l}
\Lambda \left( a;\vec I\right) \equiv \prod\limits_{r=1}^p\lambda
_{i_r}\left( a\right) \omega ^{-i_r}\,\,,\quad \quad {\rm and} \\ 
\hat \Gamma \left( \vec I\right) =\prod\limits_{r=1}^p\left( \gamma
_r^{n-1}\bar \gamma _r\right) ^{i_r}\,\,,\quad \quad {\rm while}
\end{array}
$$

\noindent (see (2.21)) 
$$
U=\prod\limits_{r=1}^p\gamma _r^{n-1}\bar \gamma _r\,\,\,, 
$$

\noindent (recall also (2.22)).\quad \quad $\quad \diamondsuit $

\smallskip\ 

\noindent This being adapted we may write: 
$$
A=\sum\limits_{\vec I\in \Gamma _p}\Lambda \left( a;\vec I\right) \hat
\Gamma \left( \vec I\right) \eqno{(A.14)}\label{A.14} 
$$

\noindent Investigation of the multisum structure of $M^q$ matrix eventuates
in rather transparent form of it, if one, (for the last time!) introduces an
overall index, or rather ''multiindex'' for all reappearing summations.

\smallskip\ 

\noindent {\bf Notation:} 
$$
\begin{array}{l}
\vec \Pi \equiv \left( k,\vec L_1,...,\vec L_q,\vec S_1,...,\vec S_q,\vec
T_1,...,\vec T_q,\vec I_1,...,\vec I_q\right) \,, \\ 
G\left( \vec \Pi \right) =\prod\limits_{r=1}^qF\left( b;k,\vec L_r,\vec
S_r,\vec T_r\right) \Lambda \left( a;\vec I_r\right) \,\,, \\ 
\hat \Omega \left( \vec \Pi ,j_1,j_2\right) =\left( \prod\limits_{r=1}^q%
\widehat{\Xi }\left( \vec L_r,\vec S_r,\vec T_r\right) \hat \Gamma \left(
\vec I_r\right) \right) U^{j_1-j_2}\omega ^{-k\left( j_1+j_2\right) }.\quad
\quad \diamondsuit
\end{array}
$$

\smallskip 

\noindent All this together taken into account leads to 
$$
M^q=\frac 1{n^2}\sum\limits_{j_1,j_2=0}^{n-1}\sum\limits_{\vec \Pi \in
\Gamma _m}G\left( \vec \Pi \right) \hat \Omega \left( \vec \Pi
,j_1,j_2\right) \eqno{(A.15)}\label{A.15} 
$$

\noindent where $\vec \Pi \in \Gamma _m\,;\quad m=3pq+1$.

\end{document}